\documentclass[a4paper,12pt]{amsart}

\usepackage[english]{babel}
\usepackage{tikz}
\usetikzlibrary{positioning,arrows,cd}

\usepackage{enumitem}

\usepackage{amsmath,amssymb,amsthm, MnSymbol}
\usepackage{lipsum}
\usepackage[utf8]{inputenc}

\usepackage{verbatim}

\newtheorem{theorem}[equation]{Theorem}
\numberwithin{equation}{section}
\newtheorem{thmx}{Theorem}

\newtheorem*{Yoneda}{Theorem}

\newtheorem{proposition}[equation]{Proposition}
\newtheorem{lemma}[equation]{Lemma}
\newtheorem{corollary}[equation]{Corollary}

\newtheorem{conjecture}[equation]{Conjecture}

\newtheorem{problem}[equation]{Problem}

\theoremstyle{definition}
\newtheorem{definition}[equation]{Definition}

\theoremstyle{remark}
				
\newtheorem{remark}[equation]{Remark}

\newcommand{\id}{\operatorname{id}}

\newcommand{\Aut}{\operatorname{Aut}}
\newcommand{\Mono}{\operatorname{Mono}}
\newcommand{\Hom}{\operatorname{Hom}}

\newcommand{\Obj}{\operatorname{Obj}}

\newcommand{\HI}{{\operatorname{HI}_\alpha}}
\newcommand{\EI}{{\operatorname{EI}_\alpha}}
\newcommand{\Rep}{{\operatorname{Rep}_\omega}}
\newcommand{\FRep}{{\operatorname{FRep}_\omega}}
\newcommand{\Epi}{\operatorname{Epi}}
\newcommand{\Ob}{\operatorname{Ob}}

\newcommand{\Bbbk}{\mathbb{K}}
\newcommand{\Z}{{\mathbb Z}}

\newcommand{\N}{{\mathbb N}}
\newcommand{\G}{\mathcal{G}}

\newcommand{\C}{\mathcal{C}}

\newcommand{\Graphs}{{{\mathcal G}raphs}}
\newcommand{\Alg}{{{\mathcal A}lg}}

\usepackage[pdfpagelabels]{hyperref}

\usepackage{colortbl}
\DefineNamedColor{named}{RoyalBlue}     {cmyk}{1,0.50,0,0}
\DefineNamedColor{named}{BrickRed}      {cmyk}{0,0.89,0.94,0.28}
\DefineNamedColor{named}{VerdeChulo}
{cmyk}{0.6,0.1,1,0}

\AtBeginDocument{%
   \def\MR#1{}
}

\begin{document}

\title[Refining Yoneda]{Hom-counting functions, combinatorial categories and related problems}

\author[A.\ Ceres]{Antonio Ceres}
\address{Departamento de \'Algebra, Geometr\'{\i}a y Topolog\'{\i}a, Universidad de M\'alaga,
Campus de Teatinos, s/n, 29071-M\'alaga, Spain}
\email{aceres@uma.es}

\author[C.\ Costoya]{Cristina Costoya}
\address{CITMAga, Departamento de Matem\'aticas,
Universidade de Santiago de Compostela, 15782-Santiago de Compostela, Spain.}
\email{cristina.costoya@usc.es}

\author[A.\ Viruel]{Antonio Viruel}
\address{Departamento de \'Algebra, Geometr\'{\i}a y Topolog\'{\i}a, Universidad de M\'alaga,
Campus de Teatinos, s/n, 29071-M\'alaga, Spain}
\email{viruel@uma.es}

\thanks{All authors were partially supported by MCIN/AEI/10.13039/501100011033 grant PID2023-149804NB-I00. The first author was also partially supported by Grant PRE2021-097724, funded by MCIN/AEI/10.13039/501100011033 and by ESF+.}

\if[]
\begin{abstract}
Classical Yoneda Lemma asserts that the isomorphism type of an object $a$ in a category $\C$ is determined by the natural type of the set-valued functor $\Hom_{\mathcal{C}}(a,\_)$. Here we show that if finiteness hypothesis are assumed to hold in $\C,$ then the isomorphism type of an object $a$ in $\C$ is determined by the integer-valued function $|\Hom_{\mathcal{C}}(a,\_)|$ on objects in $\C.$ We present applications of this result to the isomorphism problem in Group, Graph, and Ring Theory.
\end{abstract}
\fi
\begin{abstract}
Combinatorial categories satisfy a stronger form of Yoneda Lemma, namely, the isomorphism type of an object can be recovered by counting the number of homomorphisms from all other objects into it. In this work, we show that this property holds for sufficiently small categories by studying the algebra of homomorphism-counting functions. We present applications of the results to the isomorphism problem in group, graph and ring theory. 
\end{abstract}

\maketitle

\section{Introduction}\label{sec:intro}

One of the first fundamental results about category theory one may encounter is Yoneda Lemma. In its covariant form it may be formulated as follows (see \cite[p.~61]{MacLane}):

\begin{Yoneda}[Yoneda Lemma]
    Let $\mathcal{C}$ be a locally small category, and $F\colon\mathcal{C}\rightarrow \operatorname{Set}$ be a covariant functor. Then for every object $a$ in $\mathcal{C}$ there is a natural bijection
    \[\operatorname{Nat(\Hom_\mathcal{C}(a,-),F)\cong F(a)}\]
    sending a natural transformation $\eta_\ast$ to $\eta_a(\id_a)$.
\end{Yoneda}

As a consequence, natural transformations $\Hom_\mathcal{C}(a,-)\Rightarrow\Hom_\mathcal{C}(b,-)$ are induced by morphisms $b\rightarrow a$ in $\mathcal{C}$. This in particular implies that the isomorphism type of an object $a$ is completely determined by the functor $\Hom_\mathcal{C}(a,-)$. A locally finite category is called combinatorial (see \cite{Pultr}) if the integer-valued function $|\Hom_\mathcal{C}(a,-)|$ already suffices to recover $a$ up to isomorphism. 

In order to state our precise results, let us introduce the functions ${h_a,e_a,h^{a},m^{a}\colon \Ob(\mathcal{C})\rightarrow \mathbb{Z},}$ defined by  
\begin{center}
    $h_a(c):=|\Hom_\mathcal{C}(a,c)|,$ and $e_a(c):=|\Epi_\mathcal{C}(a,c)|,$
\end{center}
and dually
\begin{center}
    $h^{a}(c):=|\Hom_\mathcal{C}(c,a)|,$ and $m^a(c):=|\Mono_\mathcal{C}(c,a)|.$
\end{center} 

In this paper we study the relationship between the functions above considering the natural algebra structure on the set of integer-valued functions on $\Ob(\C)$, where addition and multiplication are defined pointwise. That is, for functions $f, g \colon \Ob(\mathcal{C}) \to \mathbb{Z}$ we define
\begin{center}
    $(f + g)(c) := f(c) + g(c),$ and $(fg)(c) := f(c)g(c),$
\end{center}
for all $c \in \Ob(\mathcal{C})$. The main result can be summarized as follows:

\begin{thmx}\label{thm:main}
Let $\C$ be a locally finite category in which morphisms admit {epi-mono factorization} and every object has a well-founded poset of subobjects and quotients. Then for any finite set of objects $\lbrace a_1,\ldots, a_n\rbrace$ in $\C$ the following are equivalent:
\begin{enumerate}[label={\rm (\roman{*})}]
    \item\label{thm:main_1}
    The objects $a_1,\dots,a_n$ are pairwise non-isomorphic.
    \item\label{thm:main_2} The functions $h_{a_1},\dots,h_{a_n}$ are linearly independent.
    \item\label{thm:main_3} The functions $e_{a_1},\dots,e_{a_n}$ are linearly independent.
    \item\label{thm:main_4} The functions $h^{a_1},\dots,h^{a_n}$ are linearly independent.
    \item\label{thm:main_5} The functions $m^{a_1},\dots,m^{a_n}$ are linearly independent.
\end{enumerate}  
\end{thmx}
\begin{proof}
    By assumption $\C$ satisfies conditions $\mathbb{M}$ and $\mathbb{E}$ (Definition~\ref{def:conditions}). Therefore, the equivalence $\ref{thm:main_2} \iff \ref{thm:main_3}$ is proven in Lemma~\ref{lem:equiv poly} while Lemma~\ref{lem:lin ind epi} shows $\ref{thm:main_2}\Rightarrow\ref{thm:main_1}$, whose converse implication being trivial.
    
    The equivalence $\ref{thm:main_1} \iff \ref{thm:main_4}\iff \ref{thm:main_5}$ follows by duality.
\end{proof}

An immediate corollary of Theorem~\ref{thm:main} is that any category satisfying the conditions above is combinatorial. 

\begin{thmx}\label{thm:main 2}
Let $\C$ be a locally finite category in which morphisms admit epi-mono factorization and every object has a well-founded poset of subobjects and quotients. Then for any given objects $a$ and $b$ in $\C$ the following are equivalent:
\begin{enumerate}[label={\rm (\roman{*})}]
    \item\label{thm:main 2_1} There is an isomorphism $a\cong b$.
    \item\label{thm:main 2_2} The functions $h_a$ and $h_b$ are equal.
    \item\label{thm:main 2_3} The functions $e_a$ and $e_b$ are equal.
    \item\label{thm:main 2_4} The functions $h^{a}$ and $h^{b}$ are equal.
    \item\label{thm:main 2_5} The functions $m^{a}$ and $m^{b}$ are equal.   
\end{enumerate}    
\end{thmx}

\begin{remark}\label{rmk:cancelacion}
As a consequence of Theorem~\ref{thm:main 2}, any category $\mathcal{C}$ satisfying its hypotheses admits both standard and power cancellation laws with respect to (co)products. That is, (co)products in $\mathcal{C}$ behave well under cancellation: if $a \otimes c \cong b \otimes c$ for some object $c$, then $a \cong b$; and similarly, if $a^{\otimes n} \cong b^{\otimes n}$ for some $n \geq 1$, then $a \cong b$ (see Section~\ref{sec:cancellation}).
\end{remark}

Finally, we show that Theorem~\ref{thm:main} admits a further extension when 
$\mathcal{C}$ has finite products or coproducts (see Theorems~\ref{thm:algebraic Yoneda epis} and \ref{thm:algebraic Yoneda monos}).

\begin{thmx}\label{thm:main 3}
Let $\mathcal{C}$ be a locally finite category in which morphisms admit epi-mono factorization and every object has a well-founded poset of subobjects and quotients. Assume that $\mathcal{C}$ has coproducts (resp.\ products). Then for any finite set of objects $\lbrace a_1,\dots,a_n\rbrace$ in $\mathcal{C}$ the following are equivalent:
\begin{enumerate}[label={\rm (\roman{*})}]
    \item\label{thm:main 3_1} Objects in $\lbrace a_1,\dots,a_n\rbrace$ are free in the semigroup $(\mathcal{C},\coprod)$ (resp.\ $(\mathcal{C},\prod)$) (see Definition \textnormal{\ref{def:monoid}}).
    \item\label{thm:main 3_2} The functions $h_{a_1},\dots,h_{a_n}$ (resp.\ $h^{a_1},\dots,h^{a_n}$) are algebraically independent.
\end{enumerate}
\end{thmx}

\subsection*{Organization of the paper}

In Section~\ref{sec:basic}, we review several elementary notions and results from category theory that will be used frequently throughout the paper, including the definitions of subobjects, quotients, and (co-)Hopfian objects. This section is included for completeness and to make the paper more self-contained; however, readers already familiar with the basic concepts of category theory may safely skip it.

Section~\ref{sec:fact_morphisms} is devoted to the factorization of morphisms and includes the key Lemma~\ref{lem:decomposition_HI}, along with its dual, Lemma~\ref{lem:decomposition_Rep}. In Section~\ref{sec:main_result}, we establish Theorem~\ref{thm:linear Yoneda epis} and its dual, Theorem~\ref{thm:linear Yoneda monos}, which together constitute the main results of the paper.

The remainder of the paper is devoted to applications of the main results. In Section~\ref{sec:groups}, we examine the case of finite groups. Section~\ref{sec:graphs} revisits results of Lovász (see~\cite{Goodall}) concerning the isomorphism type of finite graphs. Finally, in Section~\ref{sec:cancellation}, we show that the standard and power cancellation laws hold for certain (co)cartesian categories, and we apply these results to the setting of finite-dimensional unitary algebras over finite fields.

\subsection*{Acknowledgements}
The interest in the problems addressed in this paper arose from a desire to better understand and generalize Theorem~2 in \cite{Luengo}. After a first version appeared on arXiv, the authors received generous feedback and bibliographical suggestions from Nathanael Arkor and Niall Taggart, which helped contextualize and improve the results of the paper. We are deeply grateful for their contribution.

\begin{remark}
    Theorem~\ref{thm:main} was partly proven by Yoshida in \cite[Theorem~B]{Yoshida} in the study of Burnside Rings of finite categories. Our result holds modulo subgroups of automorphisms and is true for more general categories.
    
    A similar version of Theorem~\ref{thm:main 2} was established by Reggio (see \cite[Theorems~4.3 and 4.5]{Reggio}) in terms of polyadic sets. His results work relative to families of monomorphisms and epimorphisms forming proper factorization systems. Note that any such factorization is necessarily an epi-mono factorization in our sense (see Definition~\ref{def:Epi-Mono_factorization}), and indeed, our proofs remain valid in his generality. Reggio's techniques are more involved, so we hope our simpler exposition makes the results more accessible to a wider audience. 
    
    Another perspective is provided by Fujino and Matsumoto, who developed stronger arithmetic for counting homomorphisms via inclusion-exclusion, at the cost of requiring the existence of certain pushouts (see \cite{Fujino}). 
\end{remark}

\section{Basics in Category Theory}\label{sec:basic}

We begin by recalling some basic notions from category theory that will be used later and may be helpful to readers less familiar with the subject. Most of these notions can be found in standard introductory texts, such as \cite{MacLane} or \cite{Riehl}.

\begin{definition}\label{def:mono}
Let $\mathcal{C}$ be a category. 
\begin{enumerate}[label={\rm (\roman{*})}]
\item A morphism $f\colon a\rightarrow b$ in $\mathcal{C}$ is called an \emph{epimorphism}, and denoted $f\colon a\twoheadrightarrow b$, if for every object $c$ and every pair of morphisms $g_1,g_2\colon b\to c$, if $g_1\circ f=g_2\circ f$, then $g_1=g_2.$  

\item A morphism $f\colon a\rightarrow b$ in $\mathcal{C}$ is called a \emph{monomorphism}, and denoted $f\colon a\hookrightarrow b$,  if for every object $c$ and every pair of morphisms $g_1,g_2\colon c\to a$, if $f\circ g_1=f\circ g_2$, then $g_1=g_2.$
\end{enumerate}
\end{definition}

The next result is immediate. 

\begin{proposition}\label{prop:composicion_preserva}
If $f\colon x\to y$ and $g\colon y\to z$ are epimorphisms (resp.\ monomorphisms), then their composite $g\circ f$ is also an epimorphism (resp.\ a monomorphism).  Moreover, if $g\circ f$ is an epimorphism (resp.\ a monomorphism), then so is $g$ (resp.\ $f$).
\end{proposition}

The following can be interpreted as finiteness assumptions on the objects themselves.

\begin{definition}\label{def:hopfian_object}
Let $\mathcal{C}$ be a category. 
\begin{enumerate}[label={\rm (\roman{*})}]
\item An object $a$ in $\mathcal{C}$ is said to be \emph{Hopfian} if every epimorphism $f\colon a\twoheadrightarrow a$ is an isomorphism.  
\item An object $a$ in $\mathcal{C}$ is said to be \emph{co-Hopfian} if every monomorphism $f\colon a\hookrightarrow a$ is an isomorphism.
\end{enumerate}
\end{definition} 

If $\C$ locally finite, that is, the set of morphisms between any pair of objects in $\C$ is finite, these two notions coincide.

\begin{lemma}\label{lem:finite_implies_hopf}
    Let $\mathcal{C}$ be a locally finite category. Then every object in $\mathcal{C}$ is Hopfian and co-Hopfian.
\end{lemma}
\begin{proof}

We show that every object in $\mathcal{C}$ is Hopfian. Let $a$ be an object in $\mathcal{C}$, and let $f \in \Hom_\mathcal{C}(a,a)$ be an epimorphism. We aim to show that $f$ is an isomorphism.

Consider the sequence $f, f^2, f^3, \dots$ in $\Hom_\mathcal{C}(a,a)$. Since this set is finite, the sequence must eventually repeat; that is, there exist integers $n > m$ such that $f^n = f^m$. Then
\[
f^m \circ f^{n-m} = f^n = f^m \circ \mathrm{id}_a.
\]
As $f$ is an epimorphism, so is $f^m$, and left cancellation yields $f^{n-m} = \mathrm{id}_a$. It follows that $f$ is invertible, with inverse $f^{n-m-1}$, and hence an isomorphism.

A straightforward dualization of this argument shows that every object in $\mathcal{C}$ is also co-Hopfian.
\end{proof}

\begin{definition}[{\cite[p.~19]{MacLane}, \cite[p.~12]{Riehl}}]\label{def:split-epi-mono}
A morphism $f\colon a\to b$ in $\mathcal{C}$ is called a \emph{split-epimorphism} (resp.\ \emph{split-monomorphism}), if $f$ has a right (resp.\ left) inverse $g\colon b \to a$, that is, a map such that $f\circ g = \text{id}_b$ (resp.\ $g\circ f=\text{id}_a$).   
\end{definition}

Observe that any split-epimorphism (resp.\ split-monomorphism) is, in particular, an epimorphism (resp.\ a monomorphism).

\begin{proposition}\label{prop:iso_vs_split}
Let $f\colon a\to b$ be a morphism in $\mathcal{C}$. The following are equivalent:
\begin{enumerate}[label={\rm (\roman{*})}]
\item\label{prop:iso_vs_split_1} $f$ is an isomorphism.  
\item\label{prop:iso_vs_split_2} $f$ is a monomorphism and a split-epimorphism. 
\item\label{prop:iso_vs_split_3} $f$ is an epimorphism and a split-monomorphism.
\end{enumerate}
\end{proposition}
\begin{proof}
Clearly, \ref{prop:iso_vs_split_1} implies \ref{prop:iso_vs_split_2} and \ref{prop:iso_vs_split_3}.

Assume now that \ref{prop:iso_vs_split_2} holds, and let $g$ be a right inverse to $f$, so that $f\circ g=\text{id}_b$. Then $f\circ g\circ f=\text{id}_b\circ f=f\circ\text{id}_a$. Since $f$ is a monomorphism, it follows that  $g\circ f=\text{id}_a$. That is, $f$ is an isomorphism with inverse $g$, and \ref{prop:iso_vs_split_1} follows.

Finally, a similar argument shows that if \ref{prop:iso_vs_split_3} holds, and $g$ is a left inverse to $f$, then $f$ is an isomorphism with inverse $g$, and again \ref{prop:iso_vs_split_1} follows.
\end{proof}

The following two lemmas provide a tool to identify isomorphic Hopfian and co-Hopfian objects.

\begin{lemma}\label{lem:epi_entre_hopfian}
Let $a,b$ be two Hopfian objects in a category $\mathcal{C}$. Then $a\cong b$ if and only if there exist epimorphisms $f\colon a\twoheadrightarrow b$ and $g\colon b\twoheadrightarrow a$.     
\end{lemma}
\begin{proof}
Assume there exist two epimorphisms $f\colon a \twoheadrightarrow b$ and $g\colon b \twoheadrightarrow a$. We will show that $f$ is a split-monomorphism, and hence an isomorphism by Proposition~\ref{prop:iso_vs_split}. 

The composition $g \circ f$ is an epimorphism, and since $a$ is Hopfian, it follows that $g \circ f = \alpha$ for some $\alpha \in \Aut(a)$. This implies that $\alpha^{-1} \circ g$ is a left inverse to $f$, so $f$ is a split-monomorphism.
\end{proof}

By dualizing the preceding arguments, we obtain the following result:  
\begin{lemma}\label{lem:mono_entre_co-Hopfian}
Let $a,b$ be two co-Hopfian objects in a category $\mathcal{C}$. Then $a\cong b$ if and only if there exist monomorphisms $f\colon a\hookrightarrow b$ and $g\colon b\hookrightarrow a$.     
\end{lemma}

We now recall the abstract notions of subobjects and quotients in arbitrary categories (see~\cite[p.~126]{MacLane}).

\begin{definition}\label{def:iso_of_morphisms}
Two morphisms $i\colon a\to c,$ $j\colon b\to c$ (resp.\ $i\colon c\to a,$ $j\colon c\to b$) are called \emph{isomorphic} if there exists an isomorphism $k\colon a\to b$ such that $i=j\circ k$ (resp.\ $j=k \circ i$).
\end{definition}

\begin{definition}\label{def:subobject_and_quotients}
Let $\mathcal{C}$ be a category. 

\begin{enumerate}[label={\rm (\roman{*})}]
 
    \item A \emph{quotient} of an object $c$ in $\mathcal{C}$ is an isomorphism class of epimorphisms with source object $c$, i.e.,\ $p\colon c\twoheadrightarrow a.$
 
    \item A \emph{subobject} of an object $c$ in $\mathcal{C}$ is an isomorphism class of monomorphisms with target object $c$, i.e.,\ $i\colon a\hookrightarrow c.$
\end{enumerate}
In either case, they are said to be \emph{proper} or \emph{strict} if the corresponding morphism $p$ or $i$ is not an isomorphism.
\end{definition}

In the following, given an object $c$ of $\C$, we denote by $\mathcal{S}(c)$ (resp. $\mathcal{Q}(c)$) the set of isomorphism classes of subobjects (resp. quotients) of $c$. Whenever $\C$ is locally finite, these sets carry a natural partial order given by $d\leq d'$ if there is a monomorphism $d\hookrightarrow d'$ (resp. an epimorphism $d'\twoheadrightarrow d$). Observe that the antisymmetry property follows by Lemmas~\ref{lem:epi_entre_hopfian} and \ref{lem:mono_entre_co-Hopfian}. We write $\mathcal{S}_0(c)$ (resp. $\mathcal{Q}_0(c)$) for the subposet consisting of isomorphism classes of proper subobjects (resp. proper quotients) of $c$.

\if[]

In order to apply an induction argument on the objects of $\mathcal{C}$ we shall give the following:

\begin{definition}\label{def:size}
  Given an object $a$ in  a category $\C$ and an integer $n$, we say that the \emph{monomorphic} (resp.\ \emph{epimorphic}) \emph{size} of $a$ is at most $n$, denoted by $\mathfrak{m}(a)\leq n$ (resp.\ $\mathfrak{e}(a)\leq n$), if $a$ has at most $n$ pairwise different subobjects (resp.\ quotients). 
   
   If every object in $\C$ has finite monomorphic (resp.\ epimorphic) size, we say that $\C$ is $\mathfrak{m}$-finite (resp.\ $\mathfrak{e}$-finite).
\end{definition}

\begin{remark}\label{rmk:finitely_size_implications}
Observe that if an object $a$ in $\mathcal{C}$ has finite monomorphic (resp.\ epimorphic) size, then $a$ is co-Hopfian (resp.\ Hopfian). Indeed, if $a$ is not co-Hopfian (resp.\ not Hopfian), then there exists a monomorphism (resp.\ epimorphism) $f \colon a \to a$ that is not an isomorphism. In this case, the iterated compositions $f^n$, for $n \in \mathbb{N}$, yield infinitely many pairwise distinct proper subobjects (resp.\ quotients) of $a$, contradicting the finiteness of its monomorphic (resp.\ epimorphic) size.

However, the converse does not hold in general. For instance, in the category of abelian groups, the additive group of rational numbers $\mathbb{Q}$ is co-Hopfian \cite[p.~58]{de_la_Harpe}, yet it admits infinitely many pairwise distinct subgroups. Dually, the additive group of integers $\mathbb{Z}$ is Hopfian—since it is finitely generated and residually finite \cite[p.~55]{de_la_Harpe}—but it has infinitely many pairwise distinct quotients.
\end{remark}

\begin{lemma}\label{lem:epi respects size}
If there is an epimorphism $f\colon a\twoheadrightarrow b$, then $\mathfrak{e}(b)\leq\mathfrak{e}(a)$. Moreover, if $a$ and $b$  have finite epimorphic size, then $\mathfrak{e}(b)=\mathfrak{e}(a)$ (as ordinals) if and only if $f$ is an isomorphism.
\end{lemma}
\begin{proof}
First, observe that every quotient of $b$ yields a quotient of $a$ via composition with $f$ by Proposition~\ref{prop:composicion_preserva}. Moreover, since $f$ is an epimorphism, the assignment is injective, thus $\mathfrak{e}(b)\leq\mathfrak{e}(a)$. Finally, if $\mathfrak{e}(b)=\mathfrak{e}(a)$, then composition with $f$ is a bijection on quotients, so there must be a quotient ${q\colon b\twoheadrightarrow c}$ corresponding to $\id_a$. That yields an isomorphism $\varphi\colon c\rightarrow a$ so that $\id_a=\varphi\circ q\circ f$. Hence $f$ is a split-monomorphism and therefore an isomorphism by Proposition~\ref{prop:iso_vs_split}.

\end{proof}

An analogous lemma for monomorphic size follows mutatis mutandis from the preceding argument, and its proof is therefore omitted.

\begin{lemma}\label{lem:mono respects size}
If there is a monomorphism ${f\colon a\hookrightarrow b}$, then $\mathfrak{m}(a)\leq\mathfrak{m}(b)$. Moreover, if $a$ and $b$ have finite monomorphic size, then $\mathfrak{m}(a)=\mathfrak{m}(b)$ if and only if $f$ is an isomorphism.
\end{lemma}

\fi

\section{Factoring morphisms}\label{sec:fact_morphisms}

In this section, we study the factorization of morphisms in a category, with a focus on decomposing morphisms as compositions of maps with specific properties, such as epimorphisms and monomorphisms.

\begin{definition}\label{def:fac_unique}
Let $\mathcal{C}$ be a category. We say that a morphism factorization 
\[\begin{tikzcd}
& c\arrow[rd, "l"] & \\
    a\arrow[rr, "f"]\arrow[ru, "r"] && b
\end{tikzcd}\]
in $\mathcal{C}$ is \emph{unique up to isomorphism} if for any other given factorization 
\[\begin{tikzcd}
    a\arrow[rr, "f"]\arrow[rd, "r'"] && b \\
	& {c'}\arrow[ru, "l'"]
\end{tikzcd}\]
in $\mathcal{C}$ there exists an isomorphism $\delta\colon c\to c'$ such that the following diagram commutes
\[\begin{tikzcd}
& c\arrow[dd, "\delta"]\arrow[rd, "l"] & \\
    a\arrow[rd, "r'"]\arrow[ru, "r"] && b \\
	& {c'}\arrow[ru, "l'"]
\end{tikzcd}\]
\end{definition}

\begin{definition}[{\cite[Sect.\ 4.4]{Borceux1}}]\label{def:Epi-Mono_factorization}
We say that morphisms in a category $\mathcal{C}$ admit \emph{epi-mono factorization} if every morphism ${f\colon a\to b}$ in $\mathcal{C}$ factors as $f=i\circ p$ where $p\colon a\twoheadrightarrow c$ is an epimorphism, and $i\colon c\hookrightarrow b$ is a monomorphism, and this factorization is unique up to (necessarily unique) isomorphism. 
\end{definition}

Categories admitting epi-mono factorizations include many classical and widely studied examples. Most notably, all regular categories \cite[Chapter~2]{Borceux1}, such as the categories of sets, groups, graphs, and abelian categories.

\if[]
This notion holds in a broad class of categories, notably in regular categories \cite[Chapter~2]{Borceux1}, such as those of sets, groups, graphs, and abelian categories.
\fi

From now on, we work modulo subgroups of automorphisms. Accordingly, for a fixed category $\mathcal{C}$, let $\alpha$ (resp.\ $\omega$) denote a fixed assignment that associates to each object $a$ in $\mathcal{C}$ a subgroup of its automorphism group $\Aut_\mathcal{C}(a)$, denoted by $\alpha(a)$ (resp.\ $\omega(a)$). No further assumptions are made on the assignment $\alpha$.

Now, given objects $a, b$ in $\mathcal{C}$, there is a natural right (resp.\ left) action of $\Aut_\mathcal{C}(a)$ (resp.\ $\Aut_\mathcal{C}(b)$) on the set $\Hom_\mathcal{C}(a,b)$ via precomposition (resp.\ postcomposition). Consequently, given assignments $\alpha$ and $\omega$ as above, these actions restrict to right (resp.\ left) actions of the subgroups $\alpha(a)$ (resp.\ $\omega(b)$). We will use the following notation for the corresponding quotient sets:

\begin{align*}
    \HI(a,b) & =\Hom_\C(a,b)/\alpha(a),\\
    \EI(a, b) & =\Epi_\C(a,b)/\alpha(a),\\
    \Rep(a, b) & =\omega(b)\backslash\Hom_\C(a,b),\text{ and}\\
    \FRep(a, b) & =\omega(b)\backslash\Mono_\C(a,b),
\end{align*}
for any pair of objects $a,b$ of $\C$. Taking cardinals, we define the corresponding integer-valued functions on $\Obj(\mathcal{C})$:
\begin{align*}
    h_a^\alpha(b) & :=|\Hom_\C(a,b)/\alpha(a)|,\\
    e_a^\alpha(b) & :=|\Epi_\C(a,b)/\alpha(a)|,\\
    h^b_\omega(a) & :=|\omega(b)\backslash\Hom_\C(a,b)|,\text{ and}\\
    m^b_\omega(a) & :=|\omega(b)\backslash\Mono_\C(a,b)|.
\end{align*}

\begin{lemma}\label{lem:decomposition_HI}
Assume that morphisms in $\C$ admit {epi-mono factorization} (Definition~\textnormal{\ref{def:Epi-Mono_factorization}}). Then, for every pair of objects $c,d$ in $\C$ there is a bijection 
$$\HI(c,d)=\bigsqcup_{d'\in\mathcal{S}(d)} \EI(c,d')$$
where the disjoint union ranges over the set of subobjects of $d$ (Definition~\textnormal{\ref{def:subobject_and_quotients}}).
\end{lemma}
\begin{proof}
By hypothesis, every morphism $f\colon c\rightarrow d$ factors uniquely (up to isomorphism) as 
\[\begin{tikzcd}
    c\arrow[rr, "f"]\arrow[rd, two heads, "\overline{f}"] && d \\
	& {d'}\arrow[ru, "i", hook]
\end{tikzcd}\]
where $d'\hookrightarrow d$ is a subobject of $d$. 

Moreover, for a fixed representative  $i\colon d' \hookrightarrow d$, and two epimorphisms $\overline{f}, \overline{f}' \colon c \twoheadrightarrow d'$ such that the diagrams $f = i \circ \overline{f} = i \circ \overline{f}'$ commute,  since $i$ is a monomorphism, it follows that $\overline{f} = \overline{f}'$. Thus, each morphism $f \colon c \to d$ corresponds bijectively to a subobject of $d$ and an epimorphism from $c$ to its domain.  Since the action of $\Aut(c)$ on both $\HI$ and $\EI$ is given by precomposition with isomorphisms, this correspondence descends naturally to the associated quotient sets.
\end{proof}

When $\C$ is finite, taking cardinals yields the following:
\begin{corollary}\label{cor:decomposition_HI}
Let $\mathcal{C}$ be a locally finite category in which morphisms admit {epi-mono factorization}. For any pair of objects $c,d\in\mathcal{C}$ we have
$$h_c^\alpha(d)=\sum_{d'\in\mathcal{S}(d)}e_c^\alpha (d'),$$
where the sum ranges over the set of subobjects of $d$.
\end{corollary}

Note that, since $\C$ is locally finite, the summands above are all zero except for finitely many subobjects $d'$ of $d$. Therefore, the sum is well-defined even if the set of subobjects of $d$ is not finite.

Dual arguments to the above yield the following result:

\begin{lemma}\label{lem:decomposition_Rep}
Assume morphisms in $\C$ admit {epi-mono factorization}. Then for every pair of objects $c,d$ in $\C$ there is a bijection 
$$\Rep(c,d)=\bigsqcup_{c'\in\mathcal{Q}(c)} \FRep(c',d)$$
where the disjoint union ranges over the set of quotient objects of $c$ (Definition~\textnormal{\ref{def:subobject_and_quotients}}).
\end{lemma}

\begin{corollary}\label{cor:decomposition_Rep}
Let $\mathcal{C}$ be a locally finite category in which  morphisms admit {epi-mono factorization}. Then for any pair of objects $c,d\in\mathcal{C}$ we have
$$h^{d}_\omega(c)=\sum_{c'\in\mathcal{Q}(c)} m^d_\omega(c'),$$
where the sum ranges over the set of quotient objects of $c$.
\end{corollary}

\section{Main results on hom-counting functions and combinatorial categories}\label{sec:main_result}

This section is devoted to proving the main results of the paper. As indicated in Section~\ref{sec:intro}, these results can be thought of as a refinement of the classical Yoneda Lemma under certain smallness assumption which we now introduce.

A poset $P$ is called well-founded if every nonempty subset $S\subseteq P$ has a minimal element $m$, that means there is no element $s\in S$ such that $s<m$. Recall that one can perform induction on any well-founded poset (see \cite[Theorem~6.10]{Jech}).

\begin{definition}\label{def:conditions}
Let $\mathcal{C}$ be a locally finite category in which morphisms admit epi-mono factorizations (Definition~\ref{def:Epi-Mono_factorization}). Then:
\begin{enumerate}[label={\rm (\roman{*})}]
    \item\label{def:condition_M} We say that \emph{$\mathcal{C}$ satisfies condition $\mathbb{M}$} if, in addition, the poset of subobjects $\mathcal{S}(c)$ is well-founded for all $c\in\C$.
    \item\label{def:condition_E} We say that \emph{$\mathcal{C}$ satisfies condition $\mathbb{E}$} if, in addition, the poset of quotients $\mathcal{Q}(c)$ is well-founded for all $c\in\C$.
\end{enumerate}
\end{definition}

\begin{theorem}\label{thm:linear Yoneda epis}
Let $a_1,\dots,a_n$ be objects in a category $\C$ that satisfies condition $\mathbb{M}$. Then, the following are equivalent:
\begin{enumerate}[label={\rm (\roman{*})}]
    \item The objects $a_1,\dots,a_n$ are pairwise non-isomorphic.
    \item The functions $h_{a_1}^\alpha,\ldots,h_{a_n}^\alpha$ are linearly independent.
    \item The functions $e_{a_1}^\alpha,\dots,e_{a_n}^\alpha$ are linearly independent.
\end{enumerate}
\end{theorem}

We show the result in the following two lemmas.

\begin{lemma}\label{lem:equiv poly}
Let $ a_1,\dots,a_n$  be objects in a category $\mathcal{C}$ that satisfies the condition $\mathbb{M}$. For any given integers $\lambda_1,\dots,\lambda_n$, the following are equivalent:
\begin{enumerate}[label={\rm (\roman{*})}]
    \item\label{lem:equiv poly 1} Equation
    $\sum_{i=1}^n\lambda_i h_{a_i}^\alpha=0$
    holds.
    \item\label{lem:equiv poly 2} Equation
    $\sum_{i=1}^n\lambda_i e_{a_i}^\alpha=0$
    holds.
\end{enumerate}
\end{lemma}
\begin{proof}
By Lemma~\ref{lem:decomposition_HI}, for every object $c \in \mathcal{C}$, we may write
\[
\sum_{i=1}^n \lambda_i h_{a_i}^\alpha(c) = \sum_{i=1}^n \lambda_i \sum_{c'\in \mathcal{S}_0(c)} e_{a_i}^\alpha(c') 
= \sum_{c'\in \mathcal{S}_0(c)} \sum_{i=1}^n \lambda_i e_{a_i}^\alpha(c'),
\]
and thus \ref{lem:equiv poly 2} clearly implies \ref{lem:equiv poly 1}.

Conversely, assume that \ref{lem:equiv poly 1} holds and take an object $c\in \C$. We will show that \ref{lem:equiv poly 2} follows by induction on the poset $\mathcal{S}(c)$.

Pick any subobject $d\in \mathcal{S}(c)$ and assume \ref{lem:equiv poly 2} holds for every proper subobject $d'$ of $d$. Then, we may write
\begin{align*}
0=\sum_{i=1}^n\lambda_i h_{a_i}^\alpha(d) & =\sum_{i=1}^n\lambda_i\sum_{d'\in \mathcal{S}_0(d)} e_{a_i}^\alpha(d')+\sum_{i=1}^n\lambda_i e_{a_i}^\alpha(d)=\\
& =\sum_{d'\in \mathcal{S}_0(d)}\sum_{i=1}^n\lambda_i e_{a_i}^\alpha(d')+\sum_{i=1}^n\lambda_i e_{a_i}^\alpha(d)= \\
& =\sum_{i=1}^n\lambda_i e_{a_i}^\alpha(d)
\end{align*}
where the first summation is zero for all $d'$ by induction hypothesis. Note that the first summation is empty whenever $d$ does not have any proper subobject. Taking $d$ to be $c$ yields the desired result.
\end{proof}

\begin{lemma}\label{lem:lin ind epi}
Let $a_1,\dots,a_n$ be pairwise non-isomorphic objects of a locally finite category $\mathcal{C}$, and let $\lambda_1,\dots,\lambda_n$ be integers. Then, the equation
\[\sum_{i=1}^n\lambda_i e_{a_i}^\alpha=0\]
holds if and only if $\lambda_i=0$ for all $i=1,\dots,n$.
\end{lemma}

\begin{proof}
Consider the partial order on $\lbrace a_1,\ldots, a_n\rbrace$ given by $a_i\leq a_j$ if there is an epimorphism $a_j\twoheadrightarrow a_i$. Let $a_k$ be a maximal element with respect to this ordering. Then, by Lemma~\ref{lem:epi_entre_hopfian}, there is no epimorphism $a_i\twoheadrightarrow a_k$ for any $i\neq k$ since the objects are pairwise non-isomorphic.

Consequently, evaluating the equation above at $a_k$, we obtain
\[0 = \sum_{i=1}^n \lambda_i e_{a_i}^\alpha(a_k) = \lambda_k e_{a_k}^\alpha(a_k),\]
and hence $\lambda_k=0$.

Now we iterate this argument for second to maximal objects and repeat inductively. It follows that $\lambda_i=0$ for all $i$.
\end{proof}

As announced in Section~\ref{sec:intro}, we recover the main results from \cite{Reggio} as a consequence of Theorem~\ref{thm:linear Yoneda epis}.

\begin{theorem}\label{thm:Yoneda epis}
Let $a$ and $b$ be two objects in a category $\mathcal{C}$ that satisfies the condition $\mathbb{M}$. Then, the following are equivalent:
\begin{enumerate}[label={\rm (\roman{*})}]
    \item\label{thm:Yoneda epis_1} There is an isomorphism $a\cong b$.
    \item\label{thm:Yoneda epis_2} The functions $h_a^\alpha$ and $h_b^\alpha$ are equal.
    \item\label{thm:Yoneda epis_3} The functions $e_a^\alpha$ and $e_b^\alpha$ are equal.
\end{enumerate}
\end{theorem}
\begin{proof}
The implications $\ref{thm:Yoneda epis_2}\Rightarrow\ref{thm:Yoneda epis_1}$ and $\ref{thm:Yoneda epis_3}\Rightarrow\ref{thm:Yoneda epis_1}$ follow by Theorem~\ref{thm:linear Yoneda epis}. The converse implications are trivial. 
\end{proof}

By a completely dual argument one can obtain completely analogous results to those above. In particular, we establish the following:

\begin{theorem}\label{thm:linear Yoneda monos}
Let $a_1,\dots,a_n$ be objects in a category $\C$ that satisfies condition $\mathbb{E}$. Then, the following are equivalent:
\begin{enumerate}[label={\rm (\roman{*})}]
    \item The objects $a_1,\dots,a_n$ are pairwise non-isomorphic.
    \item The functions $h^{a_1}_\alpha,\ldots,h^{a_n}_\alpha$ are linearly independent.
    \item The functions $m^{a_1}_\alpha,\dots,m^{a_n}_\alpha$ are linearly independent.
\end{enumerate}
\end{theorem}

\begin{theorem}\label{thm:Yoneda monos}
Let $a$ and $b$ be two objects in a category $\mathcal{C}$ which satisfies condition $\mathbb{E}$. Then, the following are equivalent:
\begin{enumerate}[label={\rm (\roman{*})}]
    \item\label{thm:Yoneda monos_1} There is an isomorphism $a\cong b$.
    \item\label{thm:Yoneda monos_2} The functions $h^{a}_\omega$ and $h^{b}_\omega$ are equal.
    \item\label{thm:Yoneda monos_3} The functions $m^{a}_\omega$ and $m^{b}_\omega$ are equal.
\end{enumerate}
\end{theorem}

Taking $\alpha$ and $\omega$ to be the trivial assignments yields the statement claimed in Theorems~\ref{thm:main} and \ref{thm:main 2}.

\begin{remark}
Note that every finite poset is well-founded, therefore,  Theorem~\ref{thm:Yoneda epis} (resp. \ref{thm:Yoneda monos}) applies to locally finite categories with epi-mono factorizations in which every object has finitely many subobjects (resp. quotients). This classical result in homomorphism counting was first proved by Pultr in \cite{Pultr}.
\end{remark}

\if[] 

\begin{lemma}\label{lem:equiv poly dual}
Let $a_1,\dots,a_n$ be objects in a category $\mathcal{C}$ that satisfies condition $\mathbb{E}$. For any given integers $\lambda_1,\dots,\lambda_n$, the following are equivalent:
\begin{enumerate}[label={\rm (\roman{*})}]
    \item\label{lem:equiv poly dual 1} Equation
    $\Sigma_{i=1}^n\lambda_i h^{a_i}_\omega=0$
    holds.
    \item\label{lem:equiv poly dual 2} Equation
    $\Sigma_{i=1}^n\lambda_i m^{a_i}_\omega=0$
    holds.
\end{enumerate}
\end{lemma}

\begin{lemma}\label{lem:lin ind mono}
Let $a_1,\dots,a_n$ be pairwise non-isomorphic objects of a locally finite category $\C$ and let $\lambda_1,\dots,\lambda_n$ be integers. Then, the equation
\[\Sigma_{i=1}^n\lambda_i m^{a_i}_\omega=0\]
holds if and only if $\lambda_i=0$ for all $i=1,\dots,n$.
\end{lemma}

\fi

Finally, we consider the case of categories with (co)products.

\begin{definition}\label{def:monoid}
Let $\mathcal{C}$ be a category with coproducts (resp.\ products). The set of isomorphism classes of objects in $\mathcal{C}$, equipped with the operation induced by the categorical coproduct (resp.\ product), forms a semigroup known as the \emph{Grothendieck semigroup} of $\mathcal{C}$, denoted by $(\mathcal{C}, \coprod)$ (resp.\ $(\mathcal{C}, \prod)$).

A set of objects $\lbrace a_1,\dots,a_n\rbrace$ in $\mathcal{C}$ is said to be \emph{free} if they generate a free semigroup of $(\mathcal{C},\coprod)$ (resp.\ $(\mathcal{C},\prod)$).
\end{definition}

\begin{remark}\label{rmk:semigroup}
The \emph{Grothendieck semigroup} of a category $\mathcal{C}$ is a commutative semigroup, which becomes a monoid with unit given by the initial (resp.\ terminal) object, provided such an object exists in $\mathcal{C}$.
\end{remark}

\begin{theorem}\label{thm:algebraic Yoneda epis}
Let $a_1,\ldots,a_n$ be objects in a category $\C$ which satisfies condition $\mathbb{M}$. Assume that $\C$ has finite coproducts. Then, the following are equivalent:
\begin{enumerate}[label={\rm (\roman{*})}]
    \item The set $\lbrace a_1, \dots, a_n \rbrace$ is free in the semigroup $(\mathcal{C}, \coprod)$.   
    \item The functions $h_{a_1}, \dots, h_{a_n}$ are algebraically independent. That is, if $P \in \mathbb{Z}[X_1, \dots, X_n]$ is an integral polynomial such that
    \[
    P\big(h^{a_1}(c), \dots, h^{a_n}(c)\big) = 0 \quad \text{for all } c \in \mathcal{C},
    \]
    then $P \equiv 0$.
\end{enumerate}
\end{theorem}

\begin{proof}
By definition, the set $\lbrace a_1, \dots, a_n \rbrace$ is free if and only if distinct $n$-tuples of integers $(k_1, \dots, k_n)$ and $(l_1, \dots, l_n)$ yield non-isomorphic objects
\[
a_1^{k_1} \sqcup \cdots \sqcup a_n^{k_n} \not\cong a_1^{l_1} \sqcup \cdots \sqcup a_n^{l_n},
\]
where $a_i^{k_i}$ denotes the $k_i$-fold coproduct of $a_i$.

Now, each monomial of the form $h_{a_1}(c)^{k_1} \cdots h_{a_n}(c)^{k_n}$ can be identified with $h_{a_1^{k_1} \sqcup \cdots \sqcup a_n^{k_n}}(c)$ by the universal property of coproducts. Therefore, the claim follows from Theorem~\ref{thm:linear Yoneda monos}.
\end{proof}

A similar result holds for categories having finite products.

\begin{theorem}\label{thm:algebraic Yoneda monos}
Let $a_1,\ldots,a_n$ be objects in a category $\C$ which satisfies condition $\mathbb{E}$. Assume that $\C$ has finite products. Then, the following are equivalent:
\begin{enumerate}[label={\rm (\roman{*})}]
    \item The set $\lbrace a_1, \dots, a_n \rbrace$ is free in the semigroup $(\mathcal{C}, \prod)$.   
    \item The functions $h^{a_1}, \dots, h^{a_n}$ are algebraically independent. That is, if $P \in \mathbb{Z}[X_1, \dots, X_n]$ is an integral polynomial such that
    \[
    P\big(h^{a_1}(c), \dots, h^{a_n}(c)\big) = 0 \quad \text{for all } c \in \mathcal{C},
    \]
    then $P \equiv 0$.
\end{enumerate}
\end{theorem}

Naturally, one may ask whether these results remain valid when the equalities in Theorem~\ref{thm:Yoneda epis}.\ref{thm:Yoneda epis_2} and Theorem~\ref{thm:Yoneda monos}.\ref{thm:Yoneda monos_2} hold only on a certain collection of objects in $\mathcal{C}$. In other words, can the isomorphism type of two objects $a$ and $b$ in $\mathcal{C}$ be determined entirely by the values of the functions $h_a^\alpha$ and $h_b^\alpha$ (resp.\ $h^a_\omega$ and $h^b_\omega$) on a distinguished collection of objects of $\mathcal{C}$? This question motivates the notion of \emph{test families} of objects:

\begin{definition}\label{def:test_family}
Let $\mathcal{D}$ be a subclass of objects of $\mathcal{C}$. We say that a subclass of objects $\mathcal{T}$ of $\mathcal{C}$ is:
\begin{enumerate}[label={\rm (\roman{*})}]
    \item A \emph{left test family for $\mathcal{D}$} if for every pair of objects $a,b\in\mathcal{D}$ such that $h^{a}_\omega(c)=h^{b}_\omega(c)$ holds for every object $c\in\mathcal{T}$, then there is an isomorphism $a\cong b$.
    \item A \emph{right test family for $\mathcal{D}$} if for every pair of objects $a,b\in\mathcal{C}$ such that $h_a^\alpha(c)=h_b^\alpha(c)$ holds for every object $c\in\mathcal{T}$, then there is an isomorphism $a\cong b$.
\end{enumerate}
\end{definition}

The following sections are devoted to the study of test families in the categories of finite groups and finite graphs.

\section{On the isomorphism type of finite groups}\label{sec:groups}
The category of finite groups satisfies the hypotheses of Theorems~\ref{thm:Yoneda epis}, \ref{thm:Yoneda monos}, \ref{thm:linear Yoneda epis}, \ref{thm:linear Yoneda monos},\ref{thm:algebraic Yoneda epis}, and \ref{thm:algebraic Yoneda monos}. As a result, these theorems can be interpreted purely in group-theoretic terms. For example, Theorems~\ref{thm:linear Yoneda epis} and \ref{thm:linear Yoneda monos} recovers \cite[Theorem~2]{Luengo}.

Moreover, all of these results can be viewed as contributions to the study of the \emph{quantitative characterization of finite (simple) groups}~\cite{Shi}: the study of finite groups through numerical invariants such as the “order of elements,” the “set of element orders,” and related quantities. Naturally, the sizes of homomorphism sets are numerical invariants associated with a group, and many classical group-theoretic quantities can be expressed in terms of such sizes.

For example, the following statement reformulates the main result of \cite{order_spectrum}, showing that the collection of cyclic groups (including~$\mathbb{Z}$) forms a left test family (Definition~\ref{def:test_family}) for the class of finite simple groups.

\begin{theorem}[Vasil’ev-Grechkoseeva-Mazurov]\label{thhm:order_spectrum}
    If $G$ is a finite simple group, and $H$ is a finite group with $h^G(C)=h^H(C)$ for every (not necessarily finite) cyclic group $C$, then $H\cong G.$
\end{theorem}
\begin{proof}
First observe that 
$$|G|=h^G(\mathbb{Z})=h^H(\mathbb{Z})=|H|.$$

Second, finite cyclic groups are closed by quotients, and therefore, following the lines in the proof of $\ref{thm:Yoneda monos_2}\Rightarrow\ref{thm:Yoneda monos_3}$ in Theorem~\ref{thm:Yoneda monos}, we obtain $h^G(C)=h^H(C),$ for every finite cyclic group $C$, if and only if $m^G(C)=m^H(C),$ for every finite cyclic group $C$. Therefore
\begin{align*}
    w(G) &:=\{d\in\mathbb{Z}^+\,:\,\text{$G$ contains an element of order $d$}\}\\
    &=\{d\in\mathbb{Z}^+\,:\,m^G({\mathbb{Z}/d})>0\}\\
    &=\{d\in\mathbb{Z}^+\,:\,m^H({\mathbb{Z}/d})>0\}\\
    &=\{d\in\mathbb{Z}^+\,:\,\text{$H$ contains an element of order $d$}\}\\
    &=: w(H).
\end{align*}

In conclusion, both $G$ and $H$ have the same order, namely $|G|=|H|$, and the same spectrum, namely $w(G)=w(H).$ Since $G$ is simple, then $G\cong H$ by \cite{order_spectrum}.
\end{proof}

\begin{remark}\label{rmk:cyc not right test}
The collection of cyclic groups is not a right test family for the class of finite simple groups. To see this, observe that given a non-abelian finite simple group $G$ and an abelian group $A$, the only possible group homomorphism $f\colon G \to A$ is the trivial one. Therefore, $h_G(A) = 1$ for every abelian group $A$, and in particular, $h_G(C) = 1$ for every cyclic group $C$. In other words, for any two non-abelian finite simple groups $G_1$ and $G_2$, we have $h_{G_1}(C) = h_{G_2}(C)$ for every cyclic group $C$, while $G_1$ and $G_2$ may not be isomorphic.
\end{remark}

Returning to our results, Theorems~\ref{thm:Yoneda epis} and \ref{thm:Yoneda monos} establish that for every pair of non-isomorphic finite groups $G_1$ and $G_2$, there exists a finite group $H$ that distinguishes between them by counting group homomorphisms. While the proofs do not provide an explicit description of such a group $H$, the arguments rely on factoring morphisms through epimorphisms and monomorphisms. This suggests that it may be possible to take $H$ to be either $G_1$ or $G_2$. The following conjecture is thus natural:

\begin{conjecture}\label{conj:loca}
Let $G_1$ and $G_2$ be finite groups. The following are equivalent:
\begin{enumerate}[label={\rm (\roman{*})}]
    \item\label{conj:loca_1} There is an isomorphism $G_1\cong G_2$.
    \item\label{conj:loca_2} $h_{G_1}(G_i)=h_{G_2}(G_i),$ for $i=1,2$.
    \item\label{conj:loca_3} $h^{G_1}(G_i)=h^{G_2}(G_i),$ for $i=1,2$.
\end{enumerate}
\end{conjecture}

In terms of test families, this can be rephrased by saying that $\lbrace G_1, G_2 \rbrace$ forms a left or right test family for the pair $\lbrace G_1, G_2 \rbrace$. The following results provide supporting evidence for Conjecture~\ref{conj:loca}.

\begin{lemma}\label{lem:loca}
Let $G_1$ and $G_2$ be finite groups. The following are equivalent:
\begin{enumerate}[label={\rm (\roman{*})}]
    \item\label{lem:loca_1} There is an isomorphism $G_1\cong G_2$.
    \item\label{lem:loca_2} There exists an epimorphism $f\colon G_1\twoheadrightarrow G_2, $ and $h_{G_1}(G_1)=h_{G_2}(G_1)$.
    \item\label{lem:loca_3} There exists a monomorphism $f\colon G_1\hookrightarrow G_2,$ and $h^{G_1}(G_2)=h^{G_2}(G_2)$.
\end{enumerate}
\end{lemma}
\begin{proof}

Clearly, if \ref{lem:loca_1} holds, then so do \ref{lem:loca_2} and \ref{lem:loca_3}.

Now assume that \ref{lem:loca_3} holds. Since $f$ is a monomorphism, for any $g, k \in \Hom(G_2, G_1)$, we have $f \circ g = f \circ k$ if and only if $g = k$. In other words, post-composition with $f$ defines an injective map
\[
f_\sharp \colon \Hom(G_2, G_1) \to \Hom(G_2, G_2),
\]
which is bijective, since both sets have the same cardinality by hypothesis. In particular, the identity morphism $\id_{G_2} \in \Hom(G_2, G_2)$ lies in the image of $f_\sharp$, so there exists $g \in \Hom(G_2, G_1)$ such that $f \circ g = \id_{G_2}$. It follows that $f$ is an epimorphism, hence an isomorphism. Therefore, \ref{lem:loca_1} holds.

The implication \ref{lem:loca_2}$\Rightarrow$\ref{lem:loca_1} is completely dual.
\end{proof}

The preceding result allows us to check Conjecture~\ref{conj:loca} in specific cases. We begin by considering the situation in which one of the groups is simple.

\begin{corollary}\label{cor:conj_loca_simples}
Let $G_1$ and $G_2$ be finite groups and assume that $G_1$ is simple. If $h^{G_1}(G_i)=h^{G_2}(G_i)$ holds for $i=1,2$ (Conjecture~\textnormal{\ref{conj:loca}}.\ref{conj:loca_3}) then $G_1\cong G_2$ (Conjecture~\textnormal{\ref{conj:loca}}.\ref{conj:loca_1}).
\end{corollary}
\begin{proof}
By assumption, there is a group homomorphism $f\colon G_1\to G_2$ corresponding to $\id_{G_1}$. Since $G_1$ is simple, $f$ must be a monomorphism, and the result follows from Lemma~\ref{lem:loca}.
\end{proof}

Observe that the previous result tells us that, whenever one of the finite groups is simple, then $\lbrace G_1,G_2\rbrace$ is a left test family for $\lbrace G_1,G_2\rbrace$.

We now examine the case in which one of the groups is abelian. To assist readers less familiar with the topic, we include a lemma summarizing fundamental properties of homomorphisms between abelian groups used in the arguments below.

\begin{lemma}\label{lem:basic_abelian}
 Let $A_i$, $i=1,2,3,$ be finite abelian groups.
 \begin{enumerate}[label={\rm (\roman{*})}]
    \item\label{lem:basic_abelian_1} Let $n,m$ be positive integers. Then $h_{\Z/n}(\Z/m)=\gcd(n,m).$ 
    \item\label{lem:basic_abelian_2} Finite product of abelian groups is both a categorical product and coproduct within abelian groups. Therefore 
    \begin{align*}
     h^{A_3}(A_1\times A_2)&=h^{A_3}(A_1)\cdot h^{A_3}(A_2),\\
     h_{A_1}(A_2\times A_3)&=h_{A_1}(A_2)\cdot h_{A_1}(A_3).
    \end{align*}   
    \item\label{lem:basic_abelian_3} Let $p$ be a prime, $A_1=\prod_{i=1}^s \Z/p^{n_i}$ and $A_2=\prod_{i=1}^r \Z/p.$ Then $h_{A_1}(A_2)=p^{r s}$.
    \item\label{lem:basic_abelian_4} For a given prime $p$, let $S_p(A_i)$ denote the unique $p$-sylow subgroup of $A_i$. Then $A_i=\prod_{p\text{ prime}} S_p(A_i),$
    and
    $$h_{A_1}(A_2)=\prod_{p\text{ prime}} h_{S_p(A_1)}\big(S_p(A_2)\big),$$
    and each factor $h_{S_p(A_1)}\big(S_p(A_2)\big)$ correspond to the $p$-factor in the prime factorization of $h_{A_1}(A_2)$.
    \end{enumerate}
\end{lemma}

Recall that an integer is said to be square-free if it is not divisible by any perfect square other than $1$. This means in the prime factorization of a square-free number, each prime factor appears with an exponent of $1$.

Our next result shows that Conjecture~\ref{conj:loca} holds for square-free exponent abelian groups.

\begin{corollary}\label{cor:conj_loca_square_free_abelian}
Let $G_1$ and $G_2$ be finite groups. Assume $G_2$ is abelian and its exponent is square-free. If $h_{G_1}(G_i)=h_{G_2}(G_i),$ for $i=1,2$ (Conjecture~\textnormal{\ref{conj:loca}}.\ref{conj:loca_2}), then $G_1\cong G_2$ (Conjecture \textnormal{\ref{conj:loca}}.\ref{conj:loca_1}). 
\end{corollary}
\begin{proof}
First notice that given a finite group $G$, and since $G_2$ is abelian, any map $G\to G_2$ factors uniquely through the abelianization of $G$, denoted by $G_{ab}.$ In other words $h_{G}(G_2)=h_{G_{ab}}(G_2),$   
and therefore we can apply the statements in Lemma~\ref{lem:basic_abelian} to describe $h_{G}(G_2).$

Indeed, following the notation above, we get that
$$h_{G}(G_2)=h_{G_{ab}}(G_2)=\prod_{p\text{ prime}}h_{S_p(G_{ab})}\big(S_p(G_2)\big),$$
and this is just the prime factorization of $h_{G}(G_2).$ Therefore, 
specializing the result above for $G=G_1,G_2$, and recalling that $(G_2)_{ab}=G_2,$ 
the statement of Conjecture~\textnormal{\ref{conj:loca}}.\ref{conj:loca_2} gives rise to 
\begin{equation}\label{eq:loca_square_free_abelian}
h_{S_p((G_1)_{ab})}\big(S_p(G_2)\big)=h_{S_p(G_2)}\big(S_p(G_2)\big),    
\end{equation}
for every prime $p$. Moreover, the exponent of $G_2$ is square-free, thus 
$S_p(G_2)= \prod_{i=1}^{r_p} \Z/p,$ and if we write $S_p((G_1)_{ab})=\prod_{i=1}^{s_p} \Z/p^{n_{p,i}},$ then Equation~\eqref{eq:loca_square_free_abelian} implies $r_p=s_p.$ So for every prime $p$ there exists an epimorphism
$$\psi_p\colon S_p((G_1)_{ab})=\prod_{i=1}^{r_p} \Z/p^{n_{p,i}} \twoheadrightarrow S_p(G_2)=\prod_{i=1}^{r_p} \Z/p,$$
that gives rise to an epimorphism
$$f\colon G_1\twoheadrightarrow (G_1)_{ab}=\prod_{p\text{ prime}} S_p((G_1)_{ab})\overset{\prod \psi_p}\twoheadrightarrow \prod_{p\text{ prime}} S_p(G_2) = G_2.$$
Finally, the result follows from Lemma~\ref{lem:loca}.
\end{proof}

A more optimistic conjecture is the following.

\begin{conjecture}\label{conj:muy_loca}
Let $G_1$ and $G_2$ be finite groups, and define $a_{i,j}=h_{G_i}(G_j)$ for $i,j=1,2$. The following are equivalent:
\begin{enumerate}[label={\rm (\roman{*})}]
    \item\label{conj:muy_loca_1} There is an isomorphism $G_1\cong G_2$.
    \item\label{conj:muy_loca_2} $\det\big((a_{i,j})\big)=0.$
\end{enumerate}
\end{conjecture}

Obviously, Conjecture~\ref{conj:muy_loca} implies Conjecture~\ref{conj:loca}. Although we currently lack theoretical justification for Conjecture~\ref{conj:muy_loca}, it is supported by some computational evidence.

\section{The profile of a graph}\label{sec:graphs}

Let $\Graphs$ be the category of finite undirected graphs. The problem of determining the isomorphism type of objects in $\Graphs$ by counting graph homomorphisms has already been considered in the literature under the terminology of profiles. Denote by $\G$ the skeleton of $\Graphs$. 

\begin{definition}[{\cite[Definition~5.1]{Goodall}}]\label{def:profile}
Let $\mathcal{P}$ be a subset of objects in $\G$ and consider a total order $\mathcal{P}=\lbrace P_1,\dots,P_n,\dots \rbrace$. The left $\mathcal{P}$-profile (resp.\ right $\mathcal{P}$-profile) of a graph $G$ is the sequence ${\lbrace h^G(P_i)\ |\ P_i\in\mathcal{P}\rbrace}$ (resp.\ ${\lbrace h_G(P_i)\ |\ P_i\in\mathcal{P}\rbrace}$).
\end{definition}

Since $\Graphs$ fulfills the hypothesis of Theorems~\ref{thm:Yoneda epis} and \ref{thm:Yoneda monos},
considering both $\alpha$ and $\omega$ the trivial assignation in Section~\ref{sec:fact_morphisms}, we recover \cite[Theorem~5.4]{Goodall}:  

\begin{theorem}[Lovász]\label{thm:iso_grafos}
Two finite graphs are isomorphic if and only if their respective left or right $\G$-profiles are equal.    
\end{theorem}

However, smaller profiles, i.e., smaller test families, distinguishing any pair of graphs have been studied for instance in \cite{Dvorak} and \cite{Fisk}. Let us recall some graph-theoretic terminology to illustrate this with examples. A graph $G$ is called \emph{$k$-degenerated} if every subgraph $H$ of $G$ contains some vertex of degree at most $k$. On the other hand, a graph $\Gamma$ is called \emph{bipartite} if every cycle in $\Gamma$ has even length. Then, we may reformulate \cite[Theorems~12,13 and 14]{Dvorak} as follows:

\begin{theorem}[Dvo\v rák]\label{thm:left test graphs}
The following are left test families in $\Graphs$:
\begin{enumerate}[label={\rm (\roman{*})}]
    \item The set of $2$-degenerated graphs in $\G$.
    \item The set of graphs in $\G$ with homomorphisms to a given non-bipartite graph $\Gamma$.
\end{enumerate}
\end{theorem}

Another attempt to classify graphs follows by considering graphs invariants, many of which can be described by means of profiles. A very classical example of this is the \emph{chromatic polynomial}, that is a function associating to any graph $G$ the number of $k$-colourings of $G$. It is easy to see that $k$-colourings of $G$ correspond bijectively with graph homomorphisms $G\rightarrow K_k^0$, where $K_k^l$ denotes the complete graph on $k$ vertices with $l$ loops in each vertex.  A graph is called \emph{chromatically unique} if it is determined, up to isomorphism, by its chromatic polynomial, and hence by its right $\lbrace K_k^0|\ k\in\N\rbrace$-profile. Clearly, the chromatic polynomial fails to distinguish many classes of graphs, such as trees. However, it characterizes for example cycles and complete graphs and more examples of infinite families of chromatically unique graphs are known, e.g.\ \cite{CHAO}.

Even more properties are captured by the so-called \emph{Tutte polynomial}. This is a polynomial in $\Z[x,y]$ specializing to the chromatic polynomial (see \cite[Section~2.2]{Garijo}).
The classification provided by the Tutte polynomial has been described in term of profiles in \cite[Theorem~11]{Garijo}.

\begin{theorem}[Garijo, Goodall, Ne\v set\v ril]\label{thm:Tutte as profile}
Two graphs have the same Tutte polynomial if and only if they have the same right {$\lbrace K_k^l|\ k,l\in\N\rbrace$-profiles.}
\end{theorem}

\section{Cancellation laws for commutative unitary algebras over Galois fields}\label{sec:cancellation}

Given a monoidal category $\C:=(\C, \otimes, \mathbf{1}, \alpha, \lambda, \rho)$ \cite[Section~VII.1]{MacLane}, it is natural to consider the so-called cancellation problems:

\begin{problem}[Standard cancellation Problem]\label{prob:standard_cancel}
    Given objects ${A,B,C}$ in $\C,$ does $A\otimes C\cong B\otimes C$ imply $A\cong B$?
\end{problem}

\begin{problem}[Power cancellation Problem]\label{prob:power_cancel}
    Given objects $A,B$ in $\C,$ and an integer $n>0,$ does $A^{\otimes n}\cong B^{\otimes n}$ imply $A\cong B$?
\end{problem}

If Problem~\ref{prob:standard_cancel} (resp.\ Problem~\ref{prob:power_cancel}) has always a positive answer, we say that $\C$ admits a standard (resp.\ power) cancellation law.

Recall that the following are prominent examples of (commutative) monoidal categories (see \cite[p.~163]{MacLane}).

\begin{definition}\label{def:cartesian monoidal}
Let $\mathcal{C}$ be a category with finite (co)products. The \emph{(co)cartesian monoidal product} is the monoidal product given by the internal categorical (co)product, which makes $\mathcal{C}$ into a (symmetric) monoidal category.
\end{definition}

Next we shall apply Theorem~\ref{thm:main 2} to show that categories satisfying finiteness conditions and with finite (co)products admit both standard and power cancellation laws with respect to the (co)cartesian monoidal product. This is a direct consequence of Theorems~\ref{thm:Yoneda epis} and \ref{thm:Yoneda monos}.

\begin{corollary}\label{cor:cancellation coproduct}
Let $\mathcal{C}$ be a cocartesian category satisfying condition $\mathbb{M}$ (Definition~\textnormal{\ref{def:conditions}.\ref{def:condition_M}}), and such that the set of homomorphisms between any pair of objects of $\mathcal{C}$ is non-empty. Then, for any pair of objects $a,b\in\mathcal{C}$ the following are equivalent:
\begin{enumerate}[label={\rm (\roman{*})}]
    \item\label{cor:cancellation coproduct 1} $a\cong b$.
    \item\label{cor:cancellation coproduct 2} $\bigsqcup_{i=1}^na\cong\bigsqcup_{i=1}^nb$ for some integer $n>0$.
    \item\label{cor:cancellation coproduct 3} $a\sqcup c\cong b\sqcup c$ for some object $c\in\mathcal{C}$.
\end{enumerate}
\end{corollary}

\begin{proof}
Clearly, \ref{cor:cancellation coproduct 1} implies both \ref{cor:cancellation coproduct 2} and \ref{cor:cancellation coproduct 3}.
Assume now \ref{cor:cancellation coproduct 2} holds so that, for every $d\in\mathcal{D}$, we have
$h_{\bigsqcup_{i=1}^n a}(d)=h_{\bigsqcup_{i=1}^n b}(d).$
By the universal property of coproduct we may rewrite this equality as 
$h_a(d)^n=h_b(d)^n$.
Consequently, we have $h_a=h_b$ and, by Theorem~\ref{thm:Yoneda epis}, we deduce that $a\cong b$. 

Finally, if \ref{cor:cancellation coproduct 3} holds, then we have
$h_{a\sqcup c}(d)=h_{b\sqcup c}(d)$
for every object $d\in \mathcal{D}$. Equivalently,
$h_a(d)\cdot h_c(d)=h_b(d)\cdot h_c(d).$
Since the sets of homomorphisms are non-empty, $h_c(d)\ne 0$ and we may cancel factors and obtain $h_a=h_b$, so that $a\cong b$ by Theorem~\ref{thm:Yoneda epis} as before.
\end{proof}

Again, we obtain the following dual result when considering categorical products.

\begin{corollary}\label{cor:cancellation product}
Let $\mathcal{C}$ be a cartesian category satisfying condition $\mathbb{E}$ and such that the set of homomorphisms between any pair of objects is non-empty. Then, for any pair of objects $a,b\in\mathcal{C}$ the following are equivalent:
\begin{enumerate}[label={\rm (\roman{*})}]
    \item\label{cor:cancellation product 1} $a\cong b$.
    \item\label{cor:cancellation product 2} $\Pi_{i=1}^n a\cong\Pi_{i=1}^n b$ for some integer $n>0$.
    \item\label{cor:cancellation product 3} $a \times c\cong b\times c$ for some object $c\in\mathcal{C}$.
\end{enumerate}
\end{corollary}

Let $\Bbbk$ be a finite field. The category of unitary finitely generated $\Bbbk$-algebras, denoted by $\textbf{ufg-}\Alg_\Bbbk$, endowed with the standard tensor product, is a locally finite monoidal category that does not admit a standard cancellation law in general. Indeed, in \cite{Gupta} it is shown that there exists a finitely generated unital $\Bbbk$-algebra $A$ such that $A\not \cong \Bbbk[X,Y,Z]$ in $\textbf{ufg-}\Alg_\Bbbk$, while $A\otimes \Bbbk[T]\cong \Bbbk[X,Y,Z]\otimes \Bbbk[T]=\Bbbk[X,Y,Z,T].$

The situation is significantly different when considering the category of unitary finite dimensional (as $\Bbbk$-vector space) $\Bbbk$-algebras, over a finite field $\Bbbk$. In this case, the tensor product of $\mathbb{K}$-algebras is a categorical coproduct \cite[Proposition II.4.1]{Kassel}, and we are within the hypotheses of Theorem~\ref{cor:cancellation coproduct}.

\begin{theorem}\label{thm:cancelacion_alg}
Let $\Bbbk$ be a finite field. Let $A,B$ be finite dimensional (as $\Bbbk$-vector space) commutative unitary $\Bbbk$-algebras. Then the following are equivalent:
\begin{enumerate}[label={\rm (\roman{*})}]
    \item\label{thm:cancelacion_alg_1} $A\cong B$.
    \item\label{thm:cancelacion_alg_2} For some integer $n>0$, $A^{\otimes n}\cong B^{\otimes n}$ (power cancellation law).
    \item\label{thm:cancelacion_alg_3} For some finite dimensional commutative unitary $\Bbbk$-algebra $C$, $A\otimes C\cong B\otimes C$ (standard cancellation law).
\end{enumerate}
\end{theorem}

\begin{remark}\label{rmk:cancellation_only_finite_dimension}
Observe that the conclusions of Theorem~\ref{thm:cancelacion_alg} do not hold whenever the hypothesis ``finite dimensional" is replaced by the weaker hypotheses ``finitely generated" or ``finite Krull dimension". In fact, the examples mentioned above, those constructed in \cite{Gupta} as a counter example of Zarinski Problem, are finitely generated algebras, thus they have finite Krull dimension, and do not admit cancellation. 
\end{remark}


\bibliographystyle{abbrv}
\bibliography{main}
\end{document}